\documentclass[12 pt,a4paper]{article}
\usepackage[latin1]{inputenc}
\usepackage{amsfonts}
\usepackage{amsmath}
\usepackage{amssymb}
\usepackage{graphicx}
\usepackage{hyperref}
\usepackage[all]{xy}
\usepackage[left=2.5cm, right=2.5cm, bottom=3cm]{geometry}
\usepackage{mathrsfs}

\begin{document}
\title{Condition number and matrices}
\author{Felipe Bottega Diniz}
\maketitle

\begin{quote}
\begin{center}
	{\bf Abstract} 
\end{center}
\small It is well known the concept of the condition number $\kappa(A) = \|A\|\|A^ {-1}\|$, where $A$ is a $n\times n$ real or complex matrix and the norm used is the spectral norm. Although it is very common to think in $\kappa(A)$ as ``the" condition number of $A$, the truth is that condition numbers are associated to problems, not just instance of problems. Our goal is to clarify this difference. We will introduce the general concept of condition number and apply it to the particular case of real or complex matrices. After this, we will introduce the classic condition number $\kappa(A)$ of a matrix and show some known results. 
\end{quote}	

\tableofcontents

\section{Introduction}
	\normalsize Let $X,Y$ be normed vector spaces and let $f:X\to Y$ be a function. We will consider $X$ as the input space, $Y$ the output space and $f$ the problem. For example, we can have $X = \mathbb{C}^{n\times n}$ with any norm, $Y = \mathbb{C}$ with absolute value and $f(A) = \det (A)$ for all $A\in X$. Therefore, $f$ models the problem of computing the determinant.
	
	For any input $x\in X\backslash\{0\}$ and a perturbation $\tilde{x}\in X$ of $x$, we define the absolute error of $\tilde{x}$ with respect to $x$ by $$\text{AbsError}(\tilde{x},x) = \| \tilde{x} - x \|$$ and the relative error by $$\text{RelError}(\tilde{x},x) = \frac{ \| \tilde{x} - x \| }{\|x\|}.$$  
	
	Consider similar definitions for the output space.\\
	
	In practice, when doing computations we are more interested in relative errors than absolute errors. This is because the size of the error alone is not so informative. It makes more sense to compare the size of the input absolute error with the size of the input (the same observation applies to the output error). This leads to the notion of relative error. Also, when doing computations we are interested in knowing how much the (relative) error of the output changes as the input is perturbed. This analysis tells us about the sensitivity of the problem to perturbations in the input. So if RelError$(\tilde{x},x)$ is small, we expect that RelError$(f(\tilde{x}), f(x))$ is also small. This would mean the problem is stable. Of course this not always the case, Mathematics simply likes to gives us challenges wherever we go. 
	
	One way to be sure how much $f$ is sensitive to perturbations in $x\in X$ is to consider the ratio $$\frac{ \text{RelError}(f(\tilde{x}), f(x)) }{\text{RelError}(\tilde{x},x)} $$
	for small perturbations $\tilde{x}$. The bigger is this value, more sensitive is the problem in $x$ with respect to this particular perturbation. To gain more information, we can bound the error of the perturbations by some $\delta > 0$ and consider the value $$\sup_{\text{RelError}(\tilde{x},x) \leq \delta} \ \frac{ \text{RelError}(f(\tilde{x}), f(x)) }{\text{RelError}(\tilde{x},x)}. $$
	This value tells us the worst case for perturbations with error no bigger than $\delta$. For example, if this value is $0.1$, this means the worst case we can have is $$\frac{ \text{RelError}(f(\tilde{x}), f(x)) }{\text{RelError}(\tilde{x},x)} = 0.1$$ for some perturbation $\tilde{x}$ such that $\text{RelError}(\tilde{x},x) \leq \delta$. In this case, we have that $$\text{RelError}(f(\tilde{x}), f(x)) = \frac{\text{RelError}(\tilde{x},x)}{10} \leq \frac{\delta}{10},$$
which means the error of the output is ten times smaller than the error of the input and is smaller than $\frac{\delta}{10}$. This is a valuable information, because now we know how to bound the output error knowing just a bound for the input error.\\ 
	
	Now we can make one step further. Suppose we can work with perturbations as small as we want (infinite precision). This means we can take $\delta$ above arbitrarily small. In this scenario, if the value of the ration is still 0.1, this means the sensitivity of the problem stays the same, despite we are increasing the precision. Note that this give us an important information about the sensitivity of the problem, in other words, taking the limit as $\delta \to 0$ helps us to understand the intrinsic sensitivity of the problem. This limit is the \emph{condition number} of $f$ in $x$. Formally, it is defined by $$\text{cond}^f(x) = \lim_{\delta \to 0} \sup_{\text{RelError}(\tilde{x},x) \leq \delta} \ \frac{ \text{RelError}(f(\tilde{x}), f(x)) }{\text{RelError}(\tilde{x},x)}. $$
	
	Note that the condition number is important not only because it gives us stability information, but also it's useful to obtain bounds, as showed above.\footnote{Indeed, if $\tilde{x}$ is sufficiently small, then RelError$(f(\tilde{x}), f(x)) \lessapprox \text{cond}^f(x)\cdot \text{RelError}(\tilde{x}, x)$. More precisely, RelError$(f(\tilde{x}), f(x)) \leq \text{cond}^f(x)\cdot \text{RelError}(\tilde{x}, x) + o\big(\text{RelError}(\tilde{x}, x)\big)$.} 
	
\section{Condition number of matrix problems}	
	Now we have defined what is a condition number, it should be clear that it is a little weird to talk about the condition number of a matrix alone. It would make more sense to talk about the condition number of a matrix problem. Despite that, the number $\kappa(A) = \|A\|\|A^{-1}\|$ is in fact a condition number as we will see soon, but it is worth to point out that $\kappa(A)$ is commonly referred to as ``the condition number of $A$" (without mention the matrix problem) mainly because of its role on numerical linear algebra. The fact is that $\kappa(A)$ is really pervasive in numerical linear algebra, appearing related to the condition number of many matrix problems. Originally $\kappa(A)$ was introduced and studied by Turing[1] and by Neumann and Goldstine[2].\\
	
	Before continuing, we want to generalize $\kappa(A)$. For each $r,s\in \mathbb{N}$, define $\kappa_{rs}(A) = \|A\|_{rs}\|A^{-1}\|_{sr}$, where $$\|A\|_{rs} = \sup_{x\in\mathbb{R}^n\backslash \{0\}} \frac{\|Ax\|_s}{\|x\|_r}$$ is the operator norm. In particular, for $r=s=2$ we have the classic condition number $\kappa(A)$. Although we are working in real space, all results are valid in complex space.\\
	
	Let $\Sigma = \{ A\in\mathbb{R}^{n\times n}:\ A \text{ is singular} \}$ and consider the input space $X = \mathbb{R}^{n\times n}\backslash\Sigma$ with norm $\| . \|_{rs}$. One of the most important problems in linear algebra is solving linear equations, which can be modelled by the function $f:X\times \mathbb{R}^n\to \mathbb{R}^n$, such that $f(A,b) = A^{-1}b$. Of course inverting $A$ and multiplying by $b$ is not a good way to solve this problem, but we will consider this for now. This naive approach is not good in terms of computation, but theoretically there is not wrong with it. So define the output space $Y = \mathbb{R}^{n\times n}\backslash\Sigma$ with norm $\| . \|_{sr}$ and consider the problem $g:X\to Y$, such that $g(A) = A^{-1}$.\\
	
	We will need the following lemma whose proof can be found in Bürgisser-Cucker[4].\\ 
	
	\textbf{Lemma 2.1:} Let $x\in\mathbb{R}^m$ and $y\in\mathbb{R}^n$ be such that $\|x\|_r = \|y\|_s = 1$. Then there is a matrix $B\in\mathbb{R}^{n\times m}$ such that $\|B\|_{rs} = 1$ and $Bx = y$.\\
	
	\textbf{Theorem 2.2:} cond$^g(A) = \kappa_{rs}(A)$, for all $A\in X$.\\
	
	\textbf{Proof:} Let $\tilde{A}$ be a perturbation such that $\tilde{A} = A-E$, for $E\in \mathbb{R}^{n\times n}$. Note that $$\text{RelError}(g(\tilde{A}),g(A)) = \text{RelError}((A-E)^{-1}, A^{-1}) = \frac{\|(A-E)^{-1} - A^{-1}\|_{sr}}{\|A^{-1}\|_{sr}} = $$
$$ = \frac{\|((I-EA^{-1})A)^{-1} - A^{-1}\|_{sr}}{\|A^{-1}\|_{sr}} = \frac{\|A^{-1}(I-EA^{-1})^{-1} - A^{-1}\|_{sr}}{\|A^{-1}\|_{sr}}.$$

	Remember we are considering perturbations $E$ such that $\|E\|_{rs} \to 0$ to compute the condition number. With this in mind, we can consider $\|EA^{-1}\|_{ss} < 1$. Therefore, $A^{-1}(I-EA^{-1})^{-1}-A^{-1} = A^{-1}\big(\sum_{i=0}^\infty (EA^{-1})^i\big)-A^{-1} = A^{-1}EA^{-1} + A^{-1}\sum_{i=2}^\infty (EA^{-1})^i$. From this we get $$\text{cond}^g(A) = \lim_{\delta \to 0} \sup_{\frac{\|E\|_{rs}}{\|A\|_{rs}} \leq \delta} \frac{\text{RelError}(g(\tilde{A}), g(A))}{\frac{\|E\|_{rs}}{\|A\|_{rs}}} = \lim_{\delta \to 0} \sup_{\frac{\|E\|_{rs}}{\|A\|_{rs}} \leq \delta} \frac{\|A^{-1}(I-EA^{-1})^{-1} - A^{-1}\|_{sr}}{\|A^{-1}\|_{sr}} \frac{\|A\|_{rs}}{\|E\|_{rs}} = $$
$$ = \lim_{\delta \to 0} \sup_{\frac{\|E\|_{rs}}{\|A\|_{rs}} \leq \delta} \frac{\|A^{-1}EA^{-1} + A^{-1}\sum_{i=2}^\infty (EA^{-1})^i\|_{sr}}{\|A^{-1}\|_{sr}} \frac{\|A\|_{rs}}{\|E\|_{rs}} \leq $$
$$ \leq \lim_{\delta \to 0} \sup_{\frac{\|E\|_{rs}}{\|A\|_{rs}} \leq \delta} \frac{\|A^{-1}\|^2_{sr}\|E\|_{rs} + \|A^{-1}\|_{sr}\sum_{i=2}^\infty \|E\|^i_{rs}\|A^{-1}\|^i_{sr}}{\|A^{-1}\|_{sr}} \frac{\|A\|_{rs}}{\|E\|_{rs}} = $$
$$ =  \|A\|_{rs}\cdot \lim_{\delta \to 0} \sup_{\frac{\|E\|_{rs}}{\|A\|_{rs}} \leq \delta} \Big(\|A^{-1}\|_{sr} + \sum_{i=2}^\infty \|E\|^{i-1}_{rs}\|A^{-1}\|^{i-1}_{sr}\Big) = \|A\|_{rs}\|A^{-1}\|_{sr} = \kappa_{rs}(A).$$

	To show the other inequality, let $y\in\mathbb{R}^n$ such that $\|A^{-1}y\|_r = \|A^{-1}\|_{sr}$ and $\|y\|_s = 1$. Define $x = \frac{A^{-1}y}{\|A^{-1}y\|_r}$, so $\|x\|_r = 1$. Using lemma 2.1, there is a matrix $B\in \mathbb{R}^{n\times n}$ be such that $\|B\|_{rs} = 1$ and $Bx = y$. Now define $E = \delta B$, for $\delta > 0$ small. For the perturbation $\tilde{A} = A-E$ we have that $$\kappa_{rs}(A) = \|A\|_{rs} \|A^{-1}\|_{sr} =  \|A\|_{rs} \|A^{-1}y\|_r = \|A\|_{rs} \|A^{-1}Bx\|_r = \|A\|_{rs} \left\|A^{-1}B\frac{A^{-1}y}{\|A^{-1}y\|_r}\right\|_r = $$
$$ = \|A\|_{rs} \frac{\|A^{-1}BA^{-1}y\|_r}{\|A^{-1}\|_{sr}} \leq \|A\|_{rs} \frac{\|A^{-1}BA^{-1}\|_{sr}}{\|A^{-1}\|_{sr}} = \frac{ \|A\|_{rs} \|A^{-1}EA^{-1}\|_{sr}}{\delta \|A^{-1}\|_{sr}}. $$

	It is worth noting that\footnote{Just note that $$\lim_{\|E\|_{rs}\to 0} \frac{\| A^{-1}\sum_{i=2}^\infty (EA^{-1})^i\|_{ss}}{\|E\|_{rs}} \leq \lim_{\|E\|_{rs}\to 0} \|A^{-1}\|_{sr}\sum_{i=2}^\infty \|E\|_{rs}^{i-1}\|A^{-1}\|_{sr}^i = 0. $$ }
$\|A^{-1}\sum_{i=2}^\infty (EA^{-1})^i\|_{ss} = o(\|E\|_{rs})$. Therefore, $$\text{RelError}(\tilde{A}^{-1}, A^{-1}) = \frac{\|A^{-1}EA^{-1} + A^{-1}\sum_{i=2}^\infty (EA^{-1})^i\|_{sr}}{\|A^{-1}\|_{sr}} = \frac{\|A^{-1}EA^{-1}\|_{sr}}{\|A^{-1}\|_{sr}} + o(\|E\|_{sr}),$$
which implies $$\kappa_{rs}(A) = \frac{\|A\|_{rs} \left( \text{RelError}(\tilde{A}^{-1}, A^{-1} ) - o(\|E\|_{sr}) \right) }{\delta} = $$
$$ = \frac{\|A\|_{rs} \left( \text{RelError}(\tilde{A}^{-1}, A^{-1} ) - o(\|E\|_{sr}) \right) }{\|E\|_{rs}} \leq \frac{\|A\|_{rs} \text{RelError}(\tilde{A}^{-1}, A^{-1} ) }{\|E\|_{rs}} = $$
$$ = \frac{\text{RelError}(\tilde{A}^{-1}, A^{-1} )}{\text{RelError}(\tilde{A}, A )} \leq \sup_{\text{RelError}(\tilde{A}, A ) \leq \frac{\delta}{\|A\|_{rs}}}  \frac{\text{RelError}(\tilde{A}^{-1}, A^{-1} )}{\text{RelError}(\tilde{A}, A )}. $$

	Making $\delta \to 0$, we conclude that $\kappa_{rs}(A) \leq \text{cond}^g(A)$. This completes the proof. $\hspace{0.7cm} \square$\\

	Now we will show some more examples. This will make it clear the connections between the condition number and $\kappa_{rs}(A)$ in numerical linear algebra. The reader more interested in these examples can check Trefethen[5].\\

	\textbf{Example 2.3:} Fix $A\in\mathbb{R}^{n\times m}$ non-singular and consider the problem of matrix-vector multiplication. This problem can be modelled by the function $f:\mathbb{R}^m\to \mathbb{R}^n$, such that $f(x) = Ax$. In this case, we have cond$^f(x) = \frac{\|A\|_{rs}\|x\|_r}{\|Ax\|_{s}}$. If $m = n$, then cond$^f(x) \leq \|A\|_{rs} \|A^{-1}\|_{sr}$. More precisely, cond$^f(x) = \alpha(x)\cdot\kappa_{rs}(A)$, where $\alpha(x) = \frac{\|x\|_r}{\|A^{-1}\|_{sr}\|Ax\|_s}$.\\
	
	\textbf{Example 2.4:} Fix $A\in\mathbb{R}^{n\times m}$ non-singular and consider the problem of solving the linear system $Ax = b$, where $b\in\mathbb{R}^n$ is the input. Note that $x = A^{-1}b$, therefore this problem is equivalent to make the matrix-vector multiplication $A^{-1}$ by $b$. This means the condition number of this problem is $\frac{\|A^{-1}\|_{sr} \|b\|_s}{\|A^{-1}b\|_r}$.\\
	
	\textbf{Example 2.5:} This time, fix $b\in\mathbb{R}^n$ and consider the problem of solving the linear system $Ax = b$, where $A\in\mathbb{R}^{n\times n}\backslash\Sigma$ is the input. This time the condition number is $\kappa_{rs}(A)$.\\
	
\section{Condition number and ill-posed problems}
	Consider a problem $f:X\to Y$ and an input $x\in X$. We say $f$ is \emph{ill-posed} in $x$ if cond$^f(x) = \infty$, otherwise we say $f$ is \emph{well-posed} in $x$. A very common phenomenon is to be able to write cond$^f(x)$ as the relativized inverse of the distance from $x$ to the set of ill-posed inputs. Denoting $\Sigma = \{ x\in X:\ \text{cond}^f(x) = \infty \}$, this means cond$^f(x) = \frac{\alpha}{d(x,\Sigma)}$, where $\alpha$ is a constant depending on $x$. Usually this constant is less relevant than the factor $\frac{1}{d(x,\Sigma)}$, therefore we conclude that inputs closer to $\Sigma$ are ill-conditioned, while inputs far from $\Sigma$ are well-conditioned. In the case of matrix space, we denote $d_{rs}(A,\Sigma) = \min\{\|A-E\|_{rs}:\ E\in\Sigma \}$, where $\Sigma$ is the set of singular matrices.\\
	
	The next theorem is often called \emph{Eckart-Young theorem}, although it dates back Eckart and Young[6]. More information about this history can be checked in Stewart[3].\\
		
	\textbf{Theorem 3.2:} Let $A\in\mathbb{R}^{n\times n}$ be non-singular. Then $$d_{rs}(A,\Sigma) = \frac{1}{\|A^{-1}\|_{sr}} . $$
	
	\textbf{Proof:} For any $E\in\Sigma$, there is $x\in \mathbb{R}^n\backslash\{0\}$ such that $Ex = 0$. Now note that $Ex = (A+(E-A))x = Ax + (E-A)x = 0$, so $x + A^{-1}(E-A)x = 0$. Therefore, $\|x\|_r = \|A^{-1}(E-A)x\|_r \leq \|A^{-1}\|_{sr}\|E-A\|_{rs}\|x\|_r$, which implies $\frac{1}{ \|A^{-1}\|_{sr} } \leq \|E-A\|_{rs}$. Since $E\in\Sigma$ was arbitrary, we conclude that $\frac{1}{\|A^{-1}\|_{sr}} \leq d_{rs}(A,\Sigma)$.\\
	
	To prove the other inequality, it is sufficient to find a singular matrix $E$ such that $\| A-E \|_{rs} \leq \frac{1}{\|A^{-1}\|_{sr}}$. Let $y\in\mathbb{R}^n$ be such that $\|A^{-1}y\|_r = \|A^{-1}\|_{sr}$ and $\|y\|_s = 1$. Define $x = \frac{A^{-1}y}{\|A^{-1}y\|_r}$, so $\|x\|_r = 1$. Using lemma 2.1, there is a matrix $B \in\mathbb{R}^{n\times n}$ such that $\|B\|_{rs} = 1$ and $Bx = -y$. Now define $E = \frac{B}{\|A^{-1}y\|_r}$ and note that $$(A+E)x = Ax + Ex = \frac{y}{\|A^{-1}y\|_r} + \frac{Bx}{\|A^{-1}y\|_r} =  \frac{y}{\|A^{-1}y\|_r} -  \frac{y}{\|A^{-1}y\|_r} = 0.$$ 
	This means $A+E\in\Sigma$. Now just note that $$d_{rs}(A,\Sigma) \leq \|A - (A+E) \|_{rs} = \|E\|_{rs} = \frac{\|B\|_{rs}}{\|A^{-1}y\|_r} = \frac{1}{\|A^{-1}\|_{sr}}.$$
	This completes the proof.$\hspace{11cm} \square$\\
	
	The next result is often called the \emph{condition number theorem}. In the case $A\in\Sigma$, we define $\kappa_{rs}(A) = \infty$. With this we have the following.\\ 
	
	\textbf{Corollary 3.3:} For any $A\in\mathbb{R}^{n\times n}\backslash\{0\}$, we have that $$\kappa_{rs}(A) = \frac{\|A\|_{rs}}{d_{rs}(A,\Sigma)}. $$
	
\section{Componentwise condition number}
	Suppose we are working with the input space $\mathbb{R}^n$. Given $x \in \mathbb{R}^n\backslash\{0\}$ and a perturbation $\tilde{x}$, sometimes we don't want to measure the error with $\frac{\|\tilde{x} - x\|}{\|x\|}$. Instead, we may consider more useful to work with the error in each entry, in other words, we may want to measure the error by using a function of $\tilde{x}$, such that it goes to 0 when each component error goes to 0. The most common choices are $$\max \left\{ \frac{ |\tilde{x}_1 - x_1| }{ |x_1| }, \ldots, \frac{ |\tilde{x}_n - x_n| }{ |x_n| } \right\}$$
and $$\frac{ |\tilde{x}_1 - x_1| }{ |x_1| } + \ldots + \frac{ |\tilde{x}_n - x_n| }{ |x_n| }, $$	
assuming each $x_i$ is non zero.\\

	Now let $X_1, \ldots, X_n$ be vector spaces with norms $\|.\|_{X_1}, \ldots, \|.\|_{X_n}$, respectively. Also let $x=(x_1, \ldots, x_n) \in (X_1\times\ldots\times X_n)\backslash\{0\}$ and let $\tilde{x} = (\tilde{x}_1, \ldots, \tilde{x}_n) \in X_1\times\ldots\times X_n$ be a perturbation of $x$. We define the \emph{componentwise relative error} between $x$ and $\tilde{x}$ by $$\text{RelError}_\text{Cw}(\tilde{x}, x) = \max\left\{ \frac{\| \tilde{x}_1 - x_1 \|_{X_1}}{\|x_1\|_{X_1}}, \ldots, \frac{\| \tilde{x}_n - x_n \|_{X_n}}{\|x_n\|_{X_n}} \right\} $$
or $$\text{RelError}_\text{Cw}(\tilde{x}, x) = \frac{ \|\tilde{x}_1 - x_1\| }{ \|x_1\| } + \ldots + \frac{ \|\tilde{x}_n - x_n\| }{ \|x_n\| }. $$
Again, we can define $\text{RelError}_\text{Cw}(\tilde{x}, x)$ in several ways, provided that it is a function of $\tilde{x}$ such that $\text{RelError}_\text{Cw}(\tilde{x}, x) \to 0$ as $\frac{ \|\tilde{x}_i - x_i\| }{ \|x_i\| } \to 0$, for each $i=1 \ldots n$.\\	
	
	Note that the example given at the beginning of this section gives us a particular case of componentwise relative error. In fact, just take $X_i = \mathbb{R}$ for all $i=1\ldots n$ and the absolute value as a norm in each $X_i$.\\
	
	With this new kind of error measure, we may consider more two kinds of condition number. Consider that $X_1\times \ldots \times X_m$ is the input space with norms $\|.\|_{X_i}$ for each $i=1\ldots m$ and $Y_1\times \ldots \times Y_n$ is the output space with norms $\|.\|_{Y_i}$ for each $i=1\ldots n$ (each $X_i$ and $Y_i$ is a vector space). If $f:X_1\times \ldots \times X_m \to Y_1\times \ldots \times Y_n$ is a problem of interest, we define the \emph{componentwise condition number} in $x\in (X_1\times \ldots \times X_m)\backslash\{0\}$ by $$\text{Cw}^f(x) = \lim_{\delta\to 0} \sup_{\text{RelError}_\text{Cw}(\tilde{x}, x) \leq \delta} \frac{\text{RelError}_\text{Cw}(f(\tilde{x}), f(x))}{\text{RelError}_\text{Cw}(\tilde{x}, x)}. $$  
	
	Now suppose the same input space and an output space $Y$ with norm $\|.\|_Y$. It may be the case we are in interested in measuring the output error in the normwise way, not the componentwise. In this case, if $f:X_1\times \ldots \times X_m \to Y$ is a problem of interest, we define the \emph{mixed condition number} in $x\in (X_1\times \ldots \times X_m)\backslash\{0\}$ by $$\text{M}^f(x) = \lim_{\delta\to 0} \sup_{\text{RelError}_\text{Cw}(\tilde{x}, x) \leq \delta} \frac{\text{RelError}(f(\tilde{x}), f(x))}{\text{RelError}_\text{Cw}(\tilde{x}, x)}. $$
	
	The case of condition numbers for normwise perturbations and componentwise output errors are not considered in the literature, so we won't consider this one. Also, note that the normwise error is a particular case of componentwise error when there is only one component. So the normwise condition number is a particular case of the componentwise condition number.\\
	
	Continuing the set of examples given, we may note there is left one situation for linear system solving. This is the problem when $A\in\mathbb{R}^{n\times n}\backslash\Sigma$ and $b\in\mathbb{R}^n$ are the input, and the output is $x=A^{-1}b$. To model this problem, we consider the function $f:(\mathbb{R}^{n\times n}\backslash\Sigma)\times\mathbb{R}^n \to \mathbb{R}^n$ such that $f(A,b) = A^{-1}b$. In this case, the most natural way to measure the input error is to consider the componentwise error. Therefore, the error of the input is $$ \text{RelError}_\text{Cw}((\tilde{A}, \tilde{b}),(A,b)) = \max\left\{ \frac{\| \tilde{A} - A \|_{rs}}{\|A\|_{rs}}, \frac{\| \tilde{b} - b \|_s}{\|b\|_s} \right\}. $$		
	
	For the output error, we may choose to use the normwise error. This means we are measuring the problem by the mixed condition number $$\text{M}^f(A,b) = \lim_{\delta\to 0} \sup_{\text{RelError}_\text{Cw}((\tilde{A}, \tilde{b}), (A,b)) \leq \delta} \frac{\text{RelError}(f(\tilde{A}, \tilde{b}), f(A,b))}{\text{RelError}_\text{Cw}((\tilde{A}, \tilde{b}), (A,b))}. $$
	
	The following theorem shows us the connection between $\kappa(A)_{rs}$ and $\text{M}^f(A,b)$, its proof can be found in Bürgisser-Cucker[4].\\
	
	\textbf{Theorem 4.1:} For any $A\in\mathbb{R}^{n\times n}\backslash\Sigma$ and $b\in\mathbb{R}^n\backslash\{0\}$, we have that      	$$\text{M}^f(A,b) = \kappa_{rs}(A) + \frac{\|A^{-1}\|_{sr} \|b\|_s}{\|A^{-1}b\|_r}. $$
	
	\textbf{Corollary 4.2:}  For any $A\in\mathbb{R}^{n\times n}\backslash\Sigma$ and $b\in\mathbb{R}^n\backslash\{0\}$, we have that $$\kappa_{rs}(A) \leq \text{M}^f(A,b) \leq 2\kappa_{rs}(A).$$
	
	This corollary is remarkable for it shows the mixed condition number has a little dependency on $b$. Again, the point is that we can write the condition number as a multiple of $\kappa_{rs}(A)$, where the factor multiplying is not so relevant. The biggest dependency is on $\kappa_{rs}(A)$. In particular, the most used norm is the spectral, so we will consider only $\kappa(A)$ from here. Theorem 4.1 and previous examples may lead us to the conclusion that $\kappa(A)$ is the most relevant factor when doing numerical analysis of linear systems. However this is not always the case. In fact, the case of triangular systems of linear equations provides an example in which $\kappa(A)$ turns out to be inadequate. It is long observed that triangular systems are generally solved to high accuracy in spite of being, in general, ill-conditioned. To be precise, when we say \emph{in general}, we mean a random matrix drawn from the standard Gaussian distribution. In other words, each entry of the matrix is a random real variable with distribution $N(0,1)$.\\
	
	In what follows, we will say that a triangular matrix $L\in\mathbb{R}^{n\times n}$ is a \emph{unit} lower triangular matrix if $\ell_{ii} = 1$ for $i=1\ldots n$.\\
	
	 \textbf{Theorem 4.3:} Let $L$ be a $n\times n$ random unit lower triangular matrix such that $\ell_{ij}\sim N(0,1)$ for $i>j$. Then the expected value of $\|L^{-1}\|_F$ is\footnote{$\|.\|_F$ stands for the Frobenius norm.} $$E[ \| L^{-1}\|_F ] = 2^n-1. $$
	 
	 \textbf{Proof:} Denote by $s_{ij}$ the entries of $L^{-1}$, for $1\leq i,j \leq n$. It's straightforward to verify that $s_{11} = 1$ and $s_{i1} = -\sum_{j=1}^{i-1} \ell_{ij} s_{j1}$ for $i=2\ldots n$. This implies that $s_{i1}$ is a random variable independent of $\ell_{ij}$ for all $j\geq i$.\\
	 
	 It is clear that $E[s_{11}^2] = 1$. For $i\geq 2$, we have that $$E[s_{ij}^2] = E\left[ \left( \sum_{j=1}^{i-1} \ell_{ij} s_{j1} \right)^2 \right] = E\left[ \sum_{j=1}^{i-1} \ell_{ij}^2 s_{j1}^2 + 2 \sum_{\substack {j=1\\ j < k}}^{i-1} \ell_{ij}\ell_{ik} s_{j1}s_{k1} \right] = $$
$$= \sum_{j=1}^{i-1} E\left[ \ell_{ij}^2 s_{j1}^2 \right] + 2 \sum_{\substack {j=1\\ j < k}}^{i-1} E\left[\ell_{ij}\ell_{ik} s_{j1}s_{k1} \right]. $$

	Since $s_{j1}$ is independent of $\ell_{ij}$, we have $E\left[ \ell_{ij}^2 s_{j1}^2 \right] = E\left[ \ell_{ij}^2 \right] E\left[ s_{j1}^2 \right] = E\left[ s_{j1}^2 \right]$. Also, since $\ell_{ik}$ is independent of $\ell_{ij}, s_{j1}$ and $s_{k1}$ (because $k>j$), we have $E\left[\ell_{ij}\ell_{ik} s_{j1}s_{k1} \right] = E\left[\ell_{ik} \right] E\left[\ell_{ij} s_{j1}s_{k1} \right] = 0$. Therefore, 
$$E[s_{i1}^2] =  \sum_{j=1}^{i-1} E\left[ s_{j1}^2 \right].$$

	Using $E[s_{11}] = 1$ and solving the recurrence, we have that $$ E[s_{i1}^2] = 2^{i-2},$$
for $i=2\ldots n$. Then the first column of $L^{-1}$ is such that $$E\left[ \sum_{i=1}^n s_{i1}^2 \right] = E[s_{11}^2]+\sum_{i=2}^n E[s_{i1}^2] = 1+\sum_{i=2}^n 2^{i-2} = 2^{n-1}.$$ 
	By an analogous argument one shows that the $k-$th column of $L^{-1}$ is such that $$E\left[ \sum_{i=1}^n s_{ik}^2 \right] = E[s_{kk}^2]+\sum_{i=k+1}^n E[s_{ik}^2] = 1+\sum_{i=k+1}^n 2^{i-k-1} = 2^{n-k}.$$
	Finally, just note that $$\hspace{3.5cm} E[\|L^{-1}\|_F^2] = E\left[ \sum_{i,j=1}^n s_{ij}^2 \right] = \sum_{k=1}^n 2^{n-k} = 2^n -1. \hspace{3.5cm}\square $$
	
	Now we can prove that these kind of random lower triangular systems are ill-conditioned with respect to $\kappa(L)$.\\
	
	\textbf{Corollary 4.4:} Let $L$ be a $n\times n$ random unit lower triangular matrix such that $\ell_{ij}\sim N(0,1)$ for $i>j$. Then $$E[ \kappa(L)^2 ] \geq n(2^n-1).$$
	
	\textbf{Proof:} First, note that $n \|L\| \|L^{-1}\| \geq \|L\|_F \|L^{-1}\|_F$ and that  $\|L\|_F \geq \sqrt{n}$. Therefore, $$\hspace{3cm} E[ \kappa(L)^2 ] = E[ \|L\|^2 \|L^{-1}\|^2 ] \geq E\left [ n \|L^{-1}\|_F^2 \right] = n(2^n-1). \hspace{2.3cm}\square $$ 
	
	For a general lower triangular matrix, we have the following result, whose proof can be found in Bürgisser-Cucker[4].\\
	
	\textbf{Theorem 4.5:} Let $L$ be a $n\times n$ random lower triangular matrix such that $\ell_{ij}\sim N(0,1)$ for all $i \geq j$. Then $$E[ \ln(\kappa(L)) ] \geq n\ln(2) - \ln(n) -1.$$
	
	This result still implies that general lower triangular matrices are ill-conditioned. What is interesting is that, despite this fact, we still can solve general lower triangular systems with high accuracy. Therefore the classic condition number doesn't reflect the true conditioning of this problem. The reason can be found in the error analysis itself. When doing the error analysis of triangular systems, we measure the input error componentwise, not normwise. Therefore, we should be using the componentwise condition number instead of the normwise condition number. Below we show a standard algorithm to solve lower triangular linear systems and its error analysis.\\
	
	\begin{center}
		\begin{tabular}{|l|}
			\hline
			\textbf{Algorithm 4.6}\\
			\hline\\
			\textbf{Input:} $L\in\mathbb{R}^{n\times n}\backslash\Sigma$, $b\in\mathbb{R}^n$\\
			\textbf{Preconditions:} $L$ is a lower triangular matrix\\\\
			\hline\\
			$x_1 = \frac{b_1}{\ell_{11}}$\\
			\verb|for| $i=2\ldots n$\\
			$\hspace{1cm}w = \displaystyle\sum_{j=1}^{i-1} \ell_{ij}x_j$\\
			$\hspace{1cm}x_i = \frac{b_1 - w}{\ell_{ii}}$\\
			\verb|return| $x = (x_1, \ldots, x_n)$\\\\
			\hline\\
			\textbf{Output:} $x\in\mathbb{R}^n$\\
			\textbf{Postconditions:} $Lx = b$\\\\
			\hline			
		\end{tabular} 
	\end{center}\bigskip\bigskip
	
	To make the error analysis of this algorithm, we will consider finite precision. In this case the \emph{epsilon machine} will be denoted by $\varepsilon_\text{mach}$ and the floating point representation of a vector $x\in\mathbb{R}^n$ will be denoted by $fl(x) = (fl(x_1), \ldots, fl(x_n))$. By definition, this representation satisfies $fl(x) = x(1+\delta)$, for some $\delta\in\mathbb{R}$ such that $| \delta | < \varepsilon_\text{mach}$.\\ 
	
	The following lemma is standard in numerical analysis, so we will omit its proof.\\

	\textbf{Lemma 4.7:} Let $x,y\in\mathbb{R}^n$. Then exists $\delta_1, \ldots, \delta_n \in \mathbb{R}$ such that $$fl \left( \sum_{i=1}^n x_iy_i \right) = \sum_{i=1}^n x_iy_i(1+\delta_i), $$
where $|\delta_i| \leq n\varepsilon_\text{mach}$ for $i=1\ldots n$.\\	
	
	\textbf{Theorem 4.8:} Suppose $\varepsilon_\text{mach} < \frac{1}{n}$ and denote by $fl(x)$ the computed solution of algorithm 4.6 for the inputs $L$ and $b$. Then there is a lower triangular matrix $E\in\mathbb{R}^{n\times n}$ such that $(L+E)fl(x) = b$, with $$|e_{ij}| \leq |\ell_{ij} | (n+2)\varepsilon_\text{mach}$$ for $i \geq j$.\\
	
	\textbf{Proof:} We will prove by induction on $n$. In the case $n = 1$, note that $$fl(x_1) = \frac{b_1}{\ell_{11}}(1+\delta), $$
where $| \delta | \leq \varepsilon_\text{mach}$. Now take $e_{11} = \frac{-\ell_{11}\delta}{1+\delta}$ and note that $$(L+E)fl(x) = \left( \ell_{11} + \frac{-\ell_{11}\delta}{1+\delta} \right)\frac{b_1}{\ell_{11}}(1+\delta) = b_1.$$
Finally, note that $$|e_{11}| = \frac{|\ell_{11}| |\delta|}{|1+\delta|} \leq \frac{|\ell_{11}| \varepsilon_\text{mach}}{1-|\delta|} \leq \frac{|\ell_{11}| \varepsilon_\text{mach}}{1-\varepsilon_\text{mach}} \leq |\ell_{11}|2\varepsilon_\text{mach}. $$

	Therefore, the theorem is valid for $n = 1$. For the case of $n > 1$, suppose the theorem is valid for the case $n-1$.\\

	Let $\overline{L} \in \mathbb{R}^{(n-1)\times (n-1)}$ be the matrix obtained by removing the $n$th row and the $n$th column of $L$, let $\overline{b} = (b_1, \ldots, b_{n-1})$ and $\overline{x} = (x_1, \ldots, x_{n-1})$, where $x = (x_1, \ldots, x_n)$ is the actually solution of the problem. In that case, $\overline{L}$ is lower triangular, nonsingular, and $\overline{L}\overline{x} = \overline{b}$. By the induction hypothesis, the computed solution $fl(\overline{x}) = (fl(x_1), \ldots, fl(x_n))$	 by algorithm 4.6 is such that there is a matrix $\overline{E}\in\mathbb{R}^{(n-1)\times (n-1)}$ satisfying $(\overline{L} + \overline{E})fl(\overline{x}) = \overline{b}$, with $$|\overline{e}_{ij}| \leq \frac{| \ell_{ij} | \varepsilon_\text{mach}}{1-(n-1)\varepsilon_\text{mach}} \leq |\ell_{ij}|(n + 2)\varepsilon_\text{mach}$$
for $i\geq j$.\\

	Now we use the algorithm to compute $fl(x_n)$. We have that $$fl(x_n) = fl \left( \frac{b_n - \sum_{j=1}^{n-1} \ell_{nj} fl (x_j)}{\ell_{nn}} \right) = \frac{\left(b_n - \sum_{j=1}^{n-1} \ell_{nj} fl (x_j)(1+\delta_j) \right)(1+\delta)}{\ell_{nn}(1+\delta')},$$
where $|\delta_j| \leq n\varepsilon_\text{mach}$ (by lemma 4.7) and $|\delta|, |\delta'| \leq \varepsilon_\text{mach}$ (using the definition of floating point). With this last equality we can write $(L+E)fl(x) = b$, with\\
$$\hspace{-1cm}E = \left[ \begin{array}{cccccc}
					\overline{e}_{11} & 0 & 0 & \ldots & 0 & 0\\
					\overline{e}_{21} & \overline{e}_{22} & 0 & \ldots & 0 & 0\\
					\vdots & \vdots & \vdots & \ddots & \vdots & \vdots\\
					\overline{e}_{n-1,1} & \overline{e}_{n-1,2} & \overline{e}_{n-1,3} & \ldots & \overline{e}_{n-1,n-1} & 0\\
					\ell_{n1}(\delta_1 + \delta + \delta_1\delta) & \ell_{n2}(\delta_2 + \delta + \delta_2\delta) & \ell_{n3}(\delta_3 + \delta + \delta_3\delta) & \ldots & \ell_{n,n-1}(\delta_{n-1} + \delta + \delta_{n-1}\delta) & \ell_{nn}\delta'\\
				\end{array} \right].$$\\
				
		Finally, note that $$|\ell_{nn}\delta'| \leq |\ell_{nn}|\varepsilon_\text{mach} \leq |\ell_{nn}|(n+2)\varepsilon_\text{mach}$$ and that, for each $j = 1\ldots n-1$, $$\hspace{0.5cm} |\ell_{nj}(\delta_j + \delta + \delta_j\delta)| \leq |\ell_{nj}|(|\delta_j| + |\delta| + |\delta_j| |\delta|) \leq |\ell_{nj}|(n + 1 + n\varepsilon_\text{mach})\varepsilon_\text{mach} \leq |\ell_{nj}|(n + 2)\varepsilon_\text{mach}. \hspace{0.5cm} \square$$ 
		
		The bounds of this theorem can be improved by considering a more refined version of lemma 4.7. But this is not necessary for the conclusion is the same: algorithm 4.6 is backward stable. Also, note that the error analysis is not made over the norm of $E$, instead, we measure the error componentwise. This gives us a reason of why we can solve triangular system with high accuracy, despite the fact that general triangular matrices are ill-conditioned. The point is that the classic condition number $\kappa(L)$ is just not the right way to measure the error in this situation.\\
		
		Although we considered random lower triangular matrices to conclude that the expected value of $\kappa(L)$ is large, usually what we have is a random dense matrix $A$ and a system $Ax = b$. After this we use the $QL$ factorization\footnote{We could have considered the upper triangular matrices from the very beginning of this section, the results would be the same.} to obtain a lower triangular matrix $L$. In this case is not clear what is the distribution of $L$. So let's consider this situation. Define $\Delta(\mathbb{R}^{n\times n})$ as the set of $n\times n$ lower triangular matrices and define the map $\psi:\mathbb{R}^{n \times n}\backslash\Sigma \to \Delta(\mathbb{R}^{n \times n})\backslash\Sigma$\ such that $\psi(A) = L$ is the lower triangular matrix obtained by the $QL$ factorization of $A$. Considering $\mathbb{R}^{n\times n}$ with the standard Gaussian distribution, we have that the distribution on $\Delta(\mathbb{R}^{n \times n})\backslash\Sigma$ induced by $\psi$ is the pushforward measure $\rho = \mu\circ\psi^{-1}$, where $\mu:\mathbb{R}^{n\times n} \to [0,1]$ is the probability measure defining the standard Gaussian distribution on $\mathbb{R}^{n \times n}$. The measure $\rho$ reflects the most common situation occurring in practice, and in this case $L$ is well-conditioned. We are just going to state the theorem below without its proof, but the interested reader can refer to Bürgisser-Cucker[4].\\
		
		\textbf{Theorem 4.9:} Let $A$ be a $n\times n$ random matrix such that $a_{ij} \sim N(0,1)$ for $1\leq i,j \leq n$. Then $$E[\ln (\kappa(L))] = \ln(n) + O(1),$$
where $L = \psi(A)$ and the expected value are computed with respect to the distribution $\rho$.\\

\section{Conclusions}
	Condition numbers plays a big role in numerical linear algebra	 but rarely they are studied as an object for their own sake. Instead, they are seen as tool to study numerical algorithms. Because of this view of condition numbers, it is important to clarify the difference between actually condition numbers and $\kappa(A)$. The truth is that $\kappa(A)$ is the condition number related to matrix inversion or the problem of solving linear system when the input is just the matrix. Admittedly, $\kappa(A)$ appears as an important factor in most condition numbers. In a sense $\kappa(A)$ is fundamentally connected to $A$ when we are dealing with conditioning. Although this justifies calling $\kappa(A)$ as ``the condition number of $A$", it is very important to realize this is just an abuse of language. The example of triangular systems showed how $\kappa(A)$ can be misleading, so we just can't think $\kappa(A)$ is a good measure of conditioning for all problems.\\
	
	This text follows the philosophy and notations of Bürgisser-Cucker[4]. This wonderful book is the first one where condition numbers are the principal object to be studied. It is highly recommended to those who want to get a deeper understanding of condition numbers.\\
	
\section*{References}

	[1] A.M. Turing. Rounding-off errors in matrix processes. Quarterly Journal of Mechanics and Applied Mathematics, 1:287-308, 1948.\\\\
	{}[2] J. von Neumann and H.H. Goldstine. Numerical inverting matrices of high order. Bulletin of the American Mathematical Society, 53:1021-1099, 1947.\\\\
	{}[3] G.W. Stewart. On the early history of the singular value decomposition. SIAM Review, 35(4):551-566, 1993.\\\\
	{}[4] P. Bürgisser and F. Cucker. Condition - The geometry of numerical algorithms. Grundlehren der mathematischen Wissenschaften 349, Springer, Heidelberg, 2013.	\\\\ 
	{}[5] L.N. Trefethen and D. Bau III. Numerical linear algebra. SIAM, 1997.\\\\  
	{}[6] C. Eckart and G. Young. The approximation of one matrix by another of lower rank. Psychometrika, 1(3):211-218, 1936.\\\\
	{}[7] J. W. Demmel. Applied numerical linear algebra. SIAM, 1997.\\\\	
	{}[8] N.J. Higham. Accuracy and stability of numerical algorithms - Second edition. SIAM, 2002.\\\\
	{}[9] E. Schmidt. Zur Theorie der linearen und nichtlinearen Integralgleichungen. Mathematische Annalen, 63(4):433-476, 1907.\\\\
	{}[10] H.Weyl. Das asymptotische Verteilungsgesetz der Eigenwerte linearer partieller Differentilgleichungen (mit einer Anwendung auf die Theorie der Hohlraumstrahlung). Mathematische Annalen, 71(4):441-479, 1912.\\\\
	{}[11] D. Viswanath and L.N. Trefethen. Condition Numbers of random triangular matrices. SIAM Journal on Matrix Analysis and Applications, 19:564-581, 1998.\\\\
	{}[12] J.H. Wilkinson. Rounding Errors in Algebraic Processes. Prentice Hall, New York, 1963.

\end{document}